\documentclass[conference, a4paper]{IEEEtran}
\IEEEoverridecommandlockouts

\usepackage{cite}
\usepackage{amsmath,amssymb,amsfonts}
\usepackage{algorithmic}
\usepackage{graphicx}
\usepackage{textcomp}
\usepackage{xcolor}
\usepackage{hyperref}
\usepackage{mathtools}
\usepackage{enumitem}

\newcommand{\Rank}[1]{\mathrm{rank}\left(#1\right)}
\newcommand{\Dim}[1]{\mathrm{dim}(#1)}
\newcommand{\Cauchy}[1]{\mathcal{C}(#1)}
\newcommand{\D}{\mathrm{d}}
\newcommand{\Lie}{\mathrm{L}}

\newcommand{\X}{\mathcal{X}}

\newcommand{\TX}{\mathcal{T}(\mathcal{X})}

\newcommand{\Span}[1]{\mathrm{span}\left\{#1\right\}}
\newcommand{\crktwo}{\underset{2}{\subset}}
\newcommand{\ad}[3][]{\text{ad}^{#1}_{#2}#3}
\newcommand{\zdot}{\dot{z}}
\newcommand{\zbdot}{\dot{\bar{z}}}
\newcommand{\zb}{\bar{z}}
\newcommand{\pad}[1]{\partial_{#1}}

\newcommand{\inv}{\includegraphics[width=.8em]{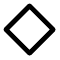}}
\newcommand{\notinv}{\includegraphics[width=.8em]{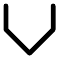}}
\newcommand{\myleftbrace}[1]{\left\{ \vphantom{\rule{0pt}{#1}} \right.}

\newtheorem{theorem}{Theorem}

\newtheorem{defn}{Definition}

\newtheorem{rem}{Remark}
\newtheorem{lem}{Lemma}

\newtheorem{exa}{Example}
\newtheorem{alg}{Algorithm}

\begin{document}

\title{ On Triangular Forms for $x$-Flat Control-Affine Systems With Two Inputs
{}
\thanks{This research was funded in whole or in part by the Austrian Science Fund (FWF) [P36473]. For open access purposes, the author has applied a CC BY public copyright license to any author accepted manuscript version arising from this submission.}
}

\author{\IEEEauthorblockN{1\textsuperscript{st} Georg Hartl}
\IEEEauthorblockA{\textit{Institute of Control Systems} \\
\textit{Johannes Kepler University} \\
Linz, Austria \\
\url{https://orcid.org/0009-0008-3797-0015}}
\and
\IEEEauthorblockN{2\textsuperscript{nd} Conrad Gstöttner}
\IEEEauthorblockA{\textit{Institute of Control Systems} \\
\textit{Johannes Kepler University}\\
Linz, Austria \\
\url{https://orcid.org/0000-0003-2107-3009}
\and
\IEEEauthorblockN{3\textsuperscript{rd} Markus Schöberl}
\IEEEauthorblockA{\textit{Institute of Control Systems} \\
\textit{Johannes Kepler University}\\
Linz, Austria \\
\url{https://orcid.org/0000-0001-5539-7015}
} } }

\maketitle

\begin{abstract}
This paper examines a broadly applicable triangular normal form for $x$-flat control-affine systems with two inputs. First, we show that this triangular form encompasses a wide range of established normal forms. Next, we prove that any $x$-flat system can be transformed into this triangular structure after a finite number of prolongations of each input. Finally, we introduce a refined algorithm for identifying candidates for $x$-flat outputs. Through illustrative examples, we demonstrate the usefulness of our results. In particular, we show that the refined algorithm exceeds the capabilities of existing methods for computing flat outputs based on triangular forms.
\end{abstract}


\section{Introduction}\label{sec1}
Differential flatness, originally introduced in the early 1990s by Fliess, Lévine, Martin, and Rouchon (see, e.g., \cite{FliessSurLesSystemes1992, FliessFlatnessDefectNonlinear1995}), has since become a central concept of nonlinear control theory by providing a powerful framework for the systematic solution of both feedforward and feedback control problems (see, e.g., \cite{FliessLieBacklundApproachEquivalence1999, GstottnerTrackingControlFlat2024}). In simple terms, a nonlinear control system 
\begin{equation}\label{eq1:nonlin_sys}
    \dot{x} = f(x,u),
\end{equation}
with the state $x=(x^1,\ldots,x^n)$ and the input $u=(u^1,\ldots,u^m)$, is flat, if there exist $m$ differentially independent functions $y^j=\varphi^j(x,u,u_{[1]},\ldots,u_{[\nu]})$, where $u_{[\nu]}$ denotes the $\nu$-th time derivative of $u$, such that $x$ and $u$ can be parameterized by $y$ and a finite number of its time derivatives. The $m$-tuple $y$ is then called a \emph{flat output} of the system. 

However, determining a flat output for a general nonlinear multi-input system remains a non-trivial task, as neither easily verifiable necessary and sufficient conditions for flatness nor computationally tractable methods for identifying flat outputs have yet been established. Exceptions include driftless two-input systems \cite{MartinFeedbackLinearizationDriftless1994}, general two-input systems that become static feedback linearizable after a one-fold prolongation of one control input \cite{NicolauFlatnessMultiInputControlAffine2017}, and those requiring a two-fold prolongation, see \cite{NicolauFlatnessTwoinputControlaffine2016,GstottnerFiniteTestLinearizability2021, GstottnerNecessarySufficientConditions2023}. 

As demonstrated, e.g., in \cite{BououdenTriangularCanonicalForm2011, GstottnerFlatTriangularForm2021, GstottnerStructurallyFlatTriangular2022, SchoberlImplicitTriangularDecomposition2014, SilveiraFlatTriangularForm2015}, structurally flat triangular forms offer an effective framework for determining flat outputs. Complete characterizations of such forms -- see, e.g.,  \cite{BououdenTriangularCanonicalForm2011, GstottnerFlatTriangularForm2021, GstottnerStructurallyFlatTriangular2022} -- yield sufficient conditions for flatness that require differentiation and algebraic operations only, thereby enabling computationally tractable methods for identifying flat outputs.

In this paper we consider the \emph{general triangular form} (GTF) introduced in \cite{GstottnerTriangularNormalForm2024} for $x$-flat control-affine two-input systems 
\begin{equation}
    \label{eq1:ai_sys}
    \dot{x} = f(x) + g_1(x)u^1 + g_2(x)u^2,
\end{equation}
where $x$-flat refers to systems that admit a flat output
\begin{equation}\label{eq1:flat_output}
    y = (\varphi^1(x), \varphi^2(x)).
\end{equation}
In \cite{GstottnerTriangularNormalForm2024} necessary and sufficient conditions were derived under which a system (\ref{eq1:ai_sys}) with a given flat output (\ref{eq1:flat_output}) can be transformed into the GTF using only state and static input transformations. We refer to systems that satisfy these conditions as \emph{static feedback equivalent} (SFE) to the GTF. Building on the GTF, we contribute the following three main results:
\begin{enumerate}
    \item We show that the GTF not only covers a wide range of established normal forms, but also includes a certain triangular form explicitly based on the so-called extended chained form. 
    \item We prove that every $x$-flat system (\ref{eq1:ai_sys}) is static feedback equivalent to the GTF after a finite number of input prolongations.
    \item We introduce a refined, distribution-based algorithm for identifying components of $x$-flat outputs. The proposed procedure extends the capabilities of some existing approaches based on triangular forms.
\end{enumerate}

Our work is structured as follows: Section \ref{sec2} introduces the notation and terminology used throughout this paper. In Section \ref{sec3}, we define differentially flat systems and present essential properties of two-input $x$-flat systems. Subsequently, Section \ref{sec4} is devoted to the general triangular form (GTF) and its relation to a triangular form based on the extended chained form. Section \ref{sec5} presents a refined algorithm for identifying components of $x$-flat outputs, followed by Section  \ref{sec6} illustrating the main contributions by several examples. Finally, the appendix provides proofs of the main theorems. \\

\section{Notation and Terminology}\label{sec2}
This work adopts tensor notation and the Einstein summation convention while omitting the index range when clear from the context. Let $\mathcal{X}$ be an $n$-dimensional smooth manifold with local coordinates $x=(x^1,\ldots,x^n)$ with its tangent space $\TX$ and let $h=(h^1,\ldots,h^m):\X \to \mathbb{R}^m$ be an $m$-tuple of smooth functions defined on $\X$. The $m\times n$ Jacobian matrix of $h$ is denoted by $\partial_xh$, and $\partial_{x^i}h^j$ represents the partial derivative of $h^j$ with respect to $x^i$. For the differentials $(\D h^1, \ldots, \D h^m)$, we write $\D h$. We denote the $k$-fold Lie derivative of a function $h^j$ along a vector field $v$ by $\Lie^k_vh^j$. 

Time derivatives are indicated by subscripts enclosed in square brackets, e.g., $h^j_{[\alpha]}$ stands for the $\alpha$-th time derivative of $h^j$, and $h_{[\alpha]} = (h^1_{[\alpha]}, \dots , h^m_{[\alpha]})$ represents the $\alpha$-th time derivative for each component of $h$. Capitalized multi-indices are used to compactly denote time derivatives of various orders across a tuple. Using the multi-index $A = (a^1, \dots, a^m)$, we obtain the concise notations $h_{[A]} = (h^1_{[a^1]}, \dots, h^m_{[a^m]})$ and $h_{[0, A]} = ( h^1_{[0,a^1]}, \dots, h^m_{[0,a^m]} )$. Here, $ h^j_{[0, a^j]}$ indicates the sequence $( h^j, \dots, h^j_{[a^j]})$ of successive time derivatives. Component-wise addition and subtraction of a multi-index with an integer are written as $A \pm c = ( a^1\pm c, \dots, a^m \pm c )$. The sum of the components of a multi-index is represented by $\#A = \sum_{j=1}^m a^j$. 

We use $[v,w]$ to indicate the Lie bracket of the vector fields $v$ and $w$. Further, the adjoint representation $\text{ad}^k_v w$ denotes the $k$-fold Lie derivative of $w$ along $v$, recursively defined by $\text{ad}^k_v w = [v, \text{ad}^{k-1}_v w]$ with $\text{ad}^0_v w = w$. Let $D_1$ and $D_2$ be two distributions. By $[D_1, D_2]$, we denote the span of the Lie brackets of  all possible pairs of the vector fields spanning $D_1$ and $D_2$, and by $[v,D_1]$ the span of the Lie brackets of the vector field $v$ and all basis vector fields of $D_1$. Given a distribution $D$, its $i$-th derived flag is defined by $D^{(i)} = D^{(i-1)} + [D^{(i-1)}, D^{(i-1)}]$, with $D^{(0)}=D$. Furthermore, its involutive closure $\bar{D}$ is the smallest involutive distribution that contains $D$ and can be determined via the derived flag. The Cauchy characteristic distribution $\Cauchy{D}$ is spanned by all vector fields $v \in D$ for which the Lie bracket with any vector field in $D$ remains in $D$, i.e., $[v,D]\subset D$.  Throughout, $\subset$ is used in the inclusive sense, meaning that equality is permitted. The notation $D_1 \underset{p}{\subset} D_2$ indicates that the corank of $D_1$ in $D_2$ is $p$. Involutive and non-involutive distributions are indicated by the symbols $\inv$ and $\notinv$, respectively.

A nonlinear system (\ref{eq1:nonlin_sys}) is said to be \emph{static feedback equivalent} (SFE) to another system $\zdot = g(z, v)$, if there exist a state transformation $z=\Phi_x(x)$ and input transformation $v=\Phi_u(x,u)$ such that the system dynamics satisfy
\vspace{-0ex}
\begin{equation*}
\left(f^i(x,u)\partial_{x^i}\Phi^j_x(x)\right)\circ\hat{\Phi}(z,v)=g^j(z,v),
\end{equation*}
\vspace{-0ex}
where $\hat{\Phi}$ denotes the inverse of $(z,v)=(\Phi_x(x),\Phi_u(x,u))$.

Finally, we assume that all functions, vector fields, and covector fields are smooth, and that all distributions and codistributions have locally constant rank. Our analysis is restricted to generic points only.

\section{ Differential Flatness }\label{sec3}

For the definition of differential flatness, we utilize a differential geometric framework as, e.g., in \cite{KolarPropertiesFlatSystems2016}, with a state-input manifold $\X \times \mathcal{U}_{[0,l_u]}$ with coordinates $(x,u_{[0,l_u]})$. Here, $l_u$ is an integer that is sufficiently large such that the time derivative of every function considered throughout this work is given by the Lie derivative along the vector field
\begin{equation}\label{eq3:vec_field_fu}
    f_u=f^i(x,u)\partial_{x^i}+\sum_{\alpha=0}^{l_u-1}u^j_{[\alpha+1]}\partial_{u^j_{[\alpha]}}.
\end{equation}
We define flatness in this setting as follows. \smallskip
\begin{defn}\label{def:flatness}
A system (\ref{eq1:nonlin_sys}) is called flat, if there exist $m$ smooth functions 
\begin{equation}\label{eq3:def_flat_output}
    \varphi^j(x,u,\ldots,u_{[\nu]}), \quad j=1,\ldots,m,
\end{equation}
defined on $\X \times \mathcal{U}_{[0,l_u]}$, such that locally the state and the input can be expressed as 
\begin{IEEEeqnarray}{rclr}
    \IEEEyesnumber\label{eq3:def_flat_para} \IEEEyessubnumber*
    x^i &\, = \,& F_x^i(\varphi_{[0,R-1]}), & \quad i=1,\ldots,n, \label{eq3:def_flat_para_x} \\
    u^j &\, = \,& F_u^j(\varphi_{[0,R]}), & \quad j=1,\ldots,m, \label{eq3:def_flat_para_u} 
\end{IEEEeqnarray}
with smooth functions $F^i_x$ and $F^j_u$, and some multi-index $R=(r_1,\ldots,r_m)$ denoting the highest orders of time derivatives occurring in (\ref{eq3:def_flat_para}). The $m$-tuple (\ref{eq3:def_flat_output}) is called a flat output of (\ref{eq1:nonlin_sys}).
\end{defn}
\smallskip
As described, e.g., in \cite{KolarPropertiesFlatSystems2016}, Definition \ref{def:flatness} implies that the differentials $\D \varphi, \D \varphi_{[1]}, \ldots, \D \varphi_{[\beta]}$ of the time derivatives of the flat output up to any order $\beta$ are linearly independent. Based upon, it can be proven that, locally, the smooth maps $F_x$ and $F_u$ along with the multi-index $R$ corresponding to a flat output (\ref{eq3:def_flat_output}) are unique. Furthermore, the map $(F_x,F_u):\mathbb{R}^{\#R+m} \to \mathbb{R}^{n+m}$ is a submersion. In, e.g., \cite{GstottnerNecessarySufficientConditions2023}, the difference between the dimension of the domain and the codomain of $(F_x,F_u)$ is called the \emph{differential difference}, represented by $d=\#R - n$.  

Focusing on $x$-flat systems, we recall an important result for two-input systems that is presented, e.g., in \cite{GstottnerLinearizationFlatTwoInput2020,GstottnerFlatSystemPossessing2023}. Consider an $x$-flat output \eqref{eq1:flat_output}. Using \eqref{eq3:vec_field_fu}, the relative degrees $K=(k_1,k_2)$ of its components are defined by
\begin{equation*}
    \Lie_{f_u}^{k_j-1}\varphi^j=\varphi^j_{[k_j-1]}(x),\quad \Lie_{f_u}^{k_j}\varphi^j=\varphi^j_{[k_j]}(x,u), \quad j=1,2.
\end{equation*}
After applying an invertible input transformation \vspace{-0.5ex}
\begin{equation}\label{eq3:input_trf}
    (\bar{u}^1, \bar{u}^2) = (\varphi^1_{[k_1]}(x,u), u^2),\vspace{-0.5ex}
\end{equation}
-- permute $u^1$ and $u^2$ if necessary -- the derivatives of the flat-output components up to the order $R$ are of the form
\begin{equation}\label{eq3:output_deriv}
    \arraycolsep=1.4pt
    \begin{array}{rclcrcl}
        \varphi^1 & = & \varphi^1(x), & \hspace{1em} & \varphi^2 & = & \varphi^2(x), \\[1ex]
         & \vdots & & & & \vdots & \\[1ex]
        \varphi^1_{[k_1-1]} & = & \varphi^1_{[k_1-1]}(x), &&  \varphi^2_{[k_2-1]} & = & \varphi^2_{[k_2-1]}(x), \\[1ex]
        \varphi^1_{[k_1]} & = & \bar{u}^1, &&  \varphi^2_{[k_2]} & = & \varphi^2_{[k_2]}(x,\bar{u}^1), \\[1ex]
        \varphi^1_{[k_1+1]} & = & \bar{u}^1_{[1]}, &&  \varphi^2_{[k_2+1]} & = & \varphi^2_{[k_2+1]}(x, \bar{u}^1_{[0,1]}), \\[1ex]
         & \vdots & & & & \vdots & \\[1ex]
        \varphi^1_{[r_1-1]} & = & \bar{u}^1_{[r_1-k_1-1]}, &&  \varphi^2_{[r_2-1]} & = & \varphi^2_{[r_2-1]}(x, \bar{u}^1_{[0,r_2-k_2-1]}), \\[1ex]
        \varphi^1_{[r_1]} & = & \bar{u}^1_{[r_1-k_1]}, &&  \varphi^2_{[r_2]} & = & \varphi^2_{[r_2]}(x, \bar{u}^1_{[0,r_2-k_2]}, \bar{u}^2) \; ,
    \end{array}
\end{equation}
and the following relations hold, see \cite{GstottnerLinearizationFlatTwoInput2020,GstottnerFlatSystemPossessing2023} for details:
\begin{subequations}
\begin{equation}\label{eq3:two_inputs_R_K_d}
    r_1-k_1 = r_2-k_2 = d ,
\end{equation}
\begin{equation}\label{eq3:two_inputs_R_K_n}
    r_1 + k_2 = n, \quad r_2 + k_1 = n, \quad n - k_1 - k_2 = d \; .
\end{equation}
\end{subequations}

\section{Triangular Forms}\label{sec4}
From now on, we consider two-input control-affine systems (\ref{eq1:ai_sys}), evolving on an $n$-dimensional state space and possessing an $x$-flat output (\ref{eq1:flat_output}). The corresponding drift and input vector fields are given by 
\begin{equation}\label{eq4:vec_fields_f_g}
    f=f^i(x)\pad{x^i}, \hspace{1em} g_1 = g^i_1(x)\pad{x^i}, \hspace{1em} g_2 = g^i_2(x)\pad{x^i}.
\end{equation}

With every $x$-flat output \eqref{eq1:flat_output}, we associate a corresponding sequence of codistributions, defined as
\begin{equation}\label{eq4:seq_codist}
    \begin{aligned}
        Q_{K-1} & = \Span{\D \varphi_{[0,K-1]}} \cap \Span{\D x}, \\
        Q_{K} & = \Span{\D \varphi_{[0,K]}} \cap \Span{\D x}, \\
        & \hspace{1ex} \vdots \\
        Q_{K+d-1} = Q_{R-1} & = \Span{\D \varphi_{[0,R-1]}} \cap \Span{\D x}.
    \end{aligned}
\end{equation}
According to Theorem 2 from \cite{GstottnerTriangularNormalForm2024}, the integrability of the codistributions \eqref{eq4:seq_codist} is a necessary and sufficient condition for a system with a given $x$-flat output $\varphi=(\varphi^1(x),\varphi^2(x))$ to be SFE to the following structurally flat triangular form, which we refer to as the \emph{general triangular form} (GTF):
\begin{equation}
    \label{eq4:triang_form_1}
    \arraycolsep=1.4pt
    \begin{array}{rl}
        \Xi_1 : \myleftbrace{6ex} & \hspace{1.4em}
        \begin{array}{rcl}
            \dot{z}^1 & = & z^2, \\
            & \hspace{0em} \vdots & \\
            \dot{z}^{k_1-1}& = & z^{k_1}, \\
        \end{array} \\
        & \hspace{2.4em} \begin{array}{rcl} \dot{z}^{k_1} & = & v^1, \end{array} \\
         \Xi_2 : \myleftbrace{6ex} & 
        \begin{array}{rcl}
            \dot{z}^{k_1+1} & = & z^{k_1+2}, \\
            &  \vdots & \\
            \dot{z}^{k_1+k_2-1} & = & z^{k_1+k_2}, \\
        \end{array} \\[5ex]
        \Xi_3 : \myleftbrace{11ex} & 
        \begin{array}{rcl}
            \dot{z}^{k_1+k_2} & = & a^{k_1+k_2}(z^1,\ldots,z^{k_1+k_2+1}) \, +  \\ 
            && \hspace{2em} b^{k_1+k_2}(z^1,\ldots,z^{k_1+k_2+1})v^1, \\
            \dot{z}^{k_1+k_2+1} & = & a^{k_1+k_2+1}(z^1,\ldots,z^{k_1+k_2+2}) \, + \\
            && \hspace{2em} b^{k_1+k_2+1}(z^1,\ldots,z^{k_1+k_2+2})v^1, \\
            &  \vdots & \\
             \dot{z}^{n-1} & = & a^{n-1}(z) + b^{n-1}(z)v^1, \\
             \dot{z}^{n} & = & v^2,
        \end{array}
    \end{array}
\end{equation}
with $b^{k_1+k_2}\neq 0$ and where, for $l=k_1+k_2,\ldots,n-1$, $\partial_{z^{l+1}}a^l\neq 0$ and/or $\partial_{z^{l+1}}b^l\neq 0$ holds. In the coordinates of the triangular form (\ref{eq4:triang_form_1}), the corresponding $x$-flat output is given by $\varphi=(z^1,z^{k_1+1})$. The drift and the input vector fields are given by
\begin{equation}\label{eq4:vec_fields_a_b}
    \begin{aligned}
        a & = z^2\pad{z^1} + \ldots + z^{k_1}\pad{z^{k_1-1}} \\ 
        & \hspace{1.3em} + z^{k_1+2}\pad{z^{k_1+1}} + \ldots +z^{k_1+k_2}\pad{z^{k_1+k_2-1}} + a^{l}\pad{z^l}, \\
        b_1 & = \pad{z^{k_1}} + b^{l}\pad{z^l}, \quad b_2 = \pad{z^n} \; ,\\
    \end{aligned}
\end{equation}
with $l=k_1+k_2,\ldots,n-1$. The GTF, as defined in (\ref{eq4:triang_form_1}), consists of two integrator chains, $\Xi_1$ and $\Xi_2$, along with a subsystem $\Xi_3$ that resembles, but is not identical to, the \emph{extended chained form} (ECF), see, e.g., \cite{KaiExtendedChainedForms2006, LiCharacterizationFlatnessExtended2013, SilveiraFlatTriangularForm2015, LiMultiinputControlaffineSystems2016}. Due to this general structure, the GTF captures several well-known normal forms, as discussed in \cite{GstottnerTriangularNormalForm2024}:
\begin{itemize}
    \item For $k_1+k_2=n$ the last subsystem reduces to \mbox{$\Xi_3: \dot{z}^n=v^2$} and we obtain the Brunovský normal form.
    \item For $\pad{z^{l+1}}b^l\neq 0$, $l=k_1+k_1,\ldots,n-1$, all functions $b^l$ can be normalized. If, in addition, the integrator chains $\Xi_1$ and $\Xi_2$ are absent (i.e., $k_1 = k_2 = 1$) and all functions $a^l$ are zero,  (\ref{eq4:triang_form_1}) reduces to the \emph{chained form} (CF), see, e.g., \cite{TilburyGoursatNormalForms1994, TilburyTrajectoryGenerationNtrailer1995, MurraySteeringNonholonomicSystems1991, MurrayNonholonomicMotionPlanning1993, MartinFeedbackLinearizationDriftless1994, LiFlatOutputsTwoinput2012}. When some functions $a^l$ are not identically zero, (\ref{eq4:triang_form_1}) corresponds to the ECF.
\end{itemize}

Furthermore, (\ref{eq4:triang_form_1}) also includes the following structurally flat triangular form, which is explicitly based on the ECF:
\vspace{-0.5ex}
\begin{equation}
    \label{eq4:triang_form_2}
    \arraycolsep=1.4pt
    \begin{array}{rl}
        \Sigma_1 : \myleftbrace{5ex} \hspace{1em} & \hspace{-0.1em}
        \begin{array}{rclcrcl}
            \dot{z}^1_{1,1} & = & z^2_{1,1}, & \hspace{1em} & \dot{z}^1_{1,2} & = & z^2_{1,2}, \\
            & \vdots & & & & \vdots & \\
            \dot{z}^{n_{1,1}}_{1,1} & = & z^1_2, & & \dot{z}^{n_{1,2}}_{1,2} & = & z^2_2, \\
        \end{array} \\[6ex]
         \Sigma_2 : \myleftbrace{9ex} \hspace{1em} & 
        \begin{array}{rcl}
           \dot{z}^1_{2} & = & z^1_{3,1},  \\
            \dot{z}^2_{2} & = & z^3_2 z^1_{3,1} + a^2_2(z_1, z^1_2, z^2_2, z^3_2), \\
            & \vdots &  \\
            \dot{z}^{n_2-1}_{2} & = & z^{n_2}_2 z^1_{3,1} + a^{n_2-1}_2(z_1, z_2), \\
            \dot{z}^{n_2}_2 & = & z^1_{3,2} + B(z_1, z_2)z^1_{3,1}, \\
        \end{array} \\[8ex]
        \Sigma_3 : \myleftbrace{9ex} \hspace{1em} & \hspace{-1em}
        \begin{array}{rclcrcl}
            \dot{z}^1_{3,1} & = & z^2_{3,1}, & & \dot{z}^1_{3,2} & = & z^2_{3,2}, \\
            & \vdots & & & & \vdots & \\
            \dot{z}^{n_3-q}_{3,1} & = & v^1, & & \dot{z}^{n_3-q}_{3,2} & = & z^{n_3-q+1}_{3,2}, \\
            &&&&& \vdots & \\
            &&& & \dot{z}^{n_3}_{3,2} & = & v^2. \\
        \end{array}
    \end{array}
\end{equation}
Subsystem $\Sigma_1$ is in Brunovský normal form and consists of two integrator chains of arbitrary lengths.\footnote{In \eqref{eq4:triang_form_2}, the first subscript of a state variable indicates whether it belongs to subsystem $\Sigma_1$, $\Sigma_2$ or $\Sigma_3$. The second subscript of state variables belonging to $\Sigma_1$ or $\Sigma_3$ refers to the integrator chain.} Subsystem $\Sigma_2$ is essentially in extended chained form, with the minor differences that the functions $a^l_2$, for $l=2,\ldots,n_2-1$, may depend on $z_1$, and that the last equation may include an additional term of the form $B(z_1, z_2)z^1_{3,1}$. Subsystem $\Sigma_3$ comprises two integrator chains whose lengths differ by an integer $q$, that is, the lengths are $n_3-q$ and $n_3$ with $n_3 \geq q$. Complete characterizations for the cases $q=0$ and $q=1$ have been established in \cite{GstottnerFlatTriangularForm2021} and \cite{GstottnerStructurallyFlatTriangular2022}, respectively. The following theorem formalizes the relationship between the triangular forms (\ref{eq4:triang_form_1}) and (\ref{eq4:triang_form_2}): 
\smallskip
\begin{theorem}\label{thm:trian_forms_1_2}
    Let system (\ref{eq1:ai_sys}) be $x$-flat and SFE to (\ref{eq4:triang_form_2}). Then it is also SFE to the GTF and the flat output $\varphi=(z^1_{1,1}, z^1_{1,2})$ of (\ref{eq4:triang_form_2}) is compatible with the GTF, in the sense that the variables $z^1_{1,1}$ and $z^1_{1,2}$ from (\ref{eq4:triang_form_2}) appear as $z^1$ and $z^{k_1+1}$ in the GTF. The corresponding relative degrees $k_1$ and $k_2$, as explicitly represented in the GTF, are then given by
    \begin{equation}\label{thm:relative_degrees}
        k_1 = n_{1,1} + n_3 - q + 1, \quad k_2 = n_{1,2} + n_3 -q + 1,
    \end{equation}
    where $n_{1,1}, n_{1,2}$ and $n_3-q$ are the lengths of the associated integrator chains of (\ref{eq4:triang_form_2}).
\end{theorem}
\smallskip
The proof of Theorem~\ref{thm:trian_forms_1_2} along with the corresponding transformation from \eqref{eq4:triang_form_2} to the GTF is provided in Appendix \ref{app:thm_triang_forms}. The following theorem highlights that the GTF is indeed a highly general structurally flat triangular form.
\smallskip
\begin{theorem}\label{thm:prolong}
    Let $y=(\varphi^1(x),\varphi^2(x))$ be an $x$-flat output of a system of the form \eqref{eq1:ai_sys}, and let $d$ be its differential difference. Then, after at most $d$-fold prolonging each input, the system becomes SFE to the GTF such that the components of $y$ appear as $(z^1,z^{k_1+1})$ in \eqref{eq4:triang_form_1}.
\end{theorem}
\smallskip
\begin{rem}\label{rem_1}
    Since a one-fold prolongation of each input of \eqref{eq1:nonlin_sys} provides a control-affine system of the form~\eqref{eq1:ai_sys}, Theorem~\ref{thm:prolong} implies that every $x$-flat system (\ref{eq1:nonlin_sys}) is SFE to the triangular form~(\ref{eq4:triang_form_1}) after at most $d+1$ prolongations of each input.
\end{rem}
\smallskip
For the proof of Theorem~\ref{thm:prolong} we refer to Appendix \ref{app:prolong}. Considering $x$-flat two-input systems it follows that \mbox{$k_1,k_2\geq1$}. Consequently, from (\ref{eq3:two_inputs_R_K_n}) it follows that 
\begin{equation}\label{eq:d_bound}
    d \leq n-2.
\end{equation}
Theorem~\ref{thm:prolong} together with Remark \ref{rem_1} and the bound \eqref{eq:d_bound} lead to the following important observation: The existence of a comprehensive characterization of the GTF -- specifically, a method to determine whether a given system is SFE to the GTF -- would enable a systematic test for $x$-flatness of systems of the form \eqref{eq1:ai_sys}. This could be done by prolonging each input $n-2$ times and then checking whether the prolonged system is SFE to the GTF.\footnote{Alternatively, one could successively one-fold prolong both inputs and test for SFE to GTF after each step, repeating this up to $n-2$ times.} Although no systematic test that checks SFE to the GTF in the general case has been developed yet, there exist algorithms that cover various subcases, for example the algorithm presented in \cite{GstottnerTriangularNormalForm2024}. An essential contribution of this work is to refine the algorithm of \cite{GstottnerTriangularNormalForm2024} in particular to cases where it is known to be insufficient. In the following section, we first review the algorithm from~\cite{GstottnerTriangularNormalForm2024}, referred to as Algorithm~\ref{alg_1}. We then analyze its limitations and introduce an improved version, referred to as Algorithm~\ref{alg_new}, which generalizes Algorithm~\ref{alg_1} and provides a straightforward method to find flat outputs for systems that are SFE to \eqref{eq4:triang_form_2} for $q\leq 1$.

\section{A refined Algorithm}\label{sec5}

Within a differential geometric framework, distribution-based methods are a well-established approach for determining flat outputs. Our work focuses on refining a distribution-based approach for finding flat outputs presented in \cite{GstottnerTriangularNormalForm2024}. Considering a system (\ref{eq1:ai_sys}) with the associated vector fields (\ref{eq4:vec_fields_f_g}) that is SFE to the GTF, we aim to construct a sequence of \mbox{\emph{feedback-invariant}} distributions\footnote{Feedback-invariant distributions are preserved under input transformations of the form $v=\Phi_u(x,u)$.}
\begin{equation}\label{eq:alg_1_sequence}
    D_1 \subset D_2 \subset \ldots \subset D_{s-1} \underset{\geq1}{\subset} D_s= \TX ,
\end{equation}
starting with $D_1=\Span{g_1,g_2}$. Based upon, we define
\begin{equation}\label{eq:alg_1_F}
    F^{\perp}=\Span{\D h^1, \ldots, \D h^c},    
\end{equation}
with $\Dim{F^\perp}=c$, functions $h^1(x), \ldots,h^c(x)$ and
\begin{itemize}
    \item $F=D_{s-1}$, if $D_{s-1}$ is involutive.
    \item $F=\Cauchy{D_{s-1}}$, if $D_{s-1}$ is non-involutive.
\end{itemize}
Consequently, candidates for flat-output components are of the form $\varphi^i(x)=\psi^i(h^1(x),\ldots,h^c(x))$ with functions \mbox{$\psi^i:\mathbb{R}^c\to\mathbb{R}$}. If one component of a flat output is known, the second component can be determined using, e.g., the method presented in Section IV of \cite{GstottnerTriangularNormalForm2024}. To verify whether a given pair of functions represents a flat output one can first determine the relative degrees $K$ and compute $R$ from Definition \ref{def:flatness} via \eqref{eq3:two_inputs_R_K_n}. It then suffices to prove that 
\begin{equation}\label{eq5:nec_suf_flat}
    \Span{\D x} \subset \Span{\D \varphi_{[0,R-1]}},
\end{equation}
which ensures that the state can be reconstructed from the flat output and its derivatives \cite{KolarContributionsDifferentialGeometric2017}.\footnote{That the input $u$ can then be expressed in terms of $\varphi_{[0,R]}$ follows from \eqref{eq1:nonlin_sys}, provided that $\Rank{\partial_u f}=m$.} 

As in \cite{GstottnerTriangularNormalForm2024}, feedback-invariant distributions (\ref{eq:alg_1_sequence}) can be computed according to the following algorithm:
\smallskip
\begin{alg}\label{alg_1}
    A sequence (\ref{eq:alg_1_sequence}) is computed by:
    \begin{enumerate}[label=\Alph*)]
        \item If $D_i$ is involutive, then $D_{i+1} = D_i + [f, D_i]$.
        \item If $D_i$ is non-involutive, then $D_{i+1} = D_i + [D_i, D_i]$.
    \end{enumerate}
\end{alg}
\smallskip

When the GTF reduces to the ECF, applying Algorithm~\ref{alg_1} yields a sequence of distributions \eqref{eq:alg_1_sequence} of the form
\begin{equation*}
    \underset{\notinv}{D_{1}} \underset{1}{\subset} \underset{\notinv}{D_{2}} \underset{1}{\subset} \ldots \underset{\notinv}{D_{n-3}} \underset{1}{\subset} \TX.
\end{equation*}
As explained, e.g., in \cite{LiFlatOutputsTwoinput2012}, $x$-flat outputs are then obtained from $\Cauchy{D_{n-3}}^{\perp}$. When the integrator chains $\Xi_1$ and/or $\Xi_2$ are present and $\Xi_3$ is essentially in ECF, we obtain the sequence 
\begin{equation*}
    \begin{aligned}
        \underset{\notinv}{D_1} \underset{1}{\subset} \underset{\notinv}{D_2} \underset{1}{\subset} & \ldots \underset{1}{\subset} \underset{\inv}{D_{n-(k_1+k_2)+1}} \crktwo \underset{\inv}{D_{n-(k_1+k_2)+2}} \crktwo \ldots \\
        &\ldots \crktwo \underset{\inv}{D_{n-\alpha}} \underset{1}{\subset} \ldots \underset{1}{\subset} \underset{\inv}{D_{n-\beta-1}} \underset{1}{\subset} D_{n-\beta} = \TX \; ,
    \end{aligned}
\end{equation*}
with $\alpha=\operatorname{max}(k_1,k_2)$ and $\beta=\operatorname{min}(k_1,k_2)$. Candidates for components of flat outputs can then be found via \mbox{$D^{\perp}_{n-\beta-1}$}. 

However, by applying Algorithm~\ref{alg_1} to systems that are SFE to the GTF it is not always possible to identify functions which are indeed flat-output components, i.e., functions which can indeed appear as $(z^1,z^{k_1+1})$ in \eqref{eq4:triang_form_1}, as the following considerations show: Consider a system of the form \eqref{eq4:triang_form_1}, where $\Xi_3$ is not in ECF and $D_1$ is involutive. Following Algorithm~\ref{alg_1}, we obtain the distributions $D_i$, \mbox{$i=1,\ldots,p+1$}, where $p$ is the smallest integer such that $D_{p+1}$ is non-involutive. Given the vector fields as defined in \eqref{eq4:vec_fields_a_b}, these distributions are of the form
\begin{equation*}
    D_{i} = \Span{b_1, \pad{n}, \ldots, \ad[i-1]{a}{b_1}, \pad{n-(i-1)} },
\end{equation*}
with the vector fields
\begin{equation}\label{eq5:vec_fields_b1}
    \ad[i]{a}{b_1} = (-1)^{i}\pad{z^{k_1 - i}} + (-1)^{i}b^{k_1 + k_2}\pad{z^{k_1 + k_2 - i}} + \ldots \; ,
\end{equation}
with $\quad b^{k_1+k_2} \neq 0 $. Consequently, the vector fields (\ref{eq5:vec_fields_b1}) contain nonzero components in $\pad{z^{k_1}},\pad{z^{k_1-1}},\ldots,\pad{z^{k_1-p}}$ and in $\pad{z^{k_1+k_2}},\pad{z^{k_1+k_2-1}},\ldots,\pad{z^{k_1+k_2-p}}$ direction. For $p \geq k_1-1$ or $p \geq k_2-1$, i.e., for $p$ too large, it can be concluded that one or even both of the differentials $(\D \varphi^1(z), \D \varphi^2(z))=(\D z^1, \D z^{k_1+1})$ might not be contained in \eqref{eq:alg_1_F}. Section~\ref{sec6} presents examples where this is the case and thus Algorithm~\ref{alg_1} does not succeed.

To motivate our refined algorithm, we consider an ideal scenario described as follows: Starting from an involutive distribution $D_1$, one constructs distributions via $ D_{i+1} = D_i + \Span{[a, \pad{z^{n - (i - 1)}}]}, i=1,\ldots,p+1$, where $p$ is the smallest integer such that $ D_{p+1} $ is non-involutive. After computing the involutive closure $\bar{D}_{p+1}$, the procedure continues by successively adding directions obtained through Lie brackets between the drift vector field $ a $ and suitable unit vector fields associated with $ \Xi_3 $, until another non-involutive distribution is reached. These two steps are repeated until all directions of $ \Xi_3 $ are spanned. Finally, the remaining directions associated with the integrator chains $\Xi_1$ and $\Xi_2$ can be added following rule A of Algorithm~\ref{alg_1}.

However, in the original coordinates $(x,u)$, identifying suitable unit vector fields associated with $\Xi_3$ is a non-trivial task. Nevertheless, we show that such vector fields can be determined under certain conditions. Therefore, we recall the key findings from the Lemmas A.5 and A.6 given in the appendix of \cite{GstottnerAnalysisControlFlat2023}, which form the basis for a refined procedure for constructing sequences of the form (\ref{eq:alg_1_sequence}).
\smallskip
\begin{lem}\label{lem_1}
    Let $D_0 \crktwo D_1 \crktwo D_2$ be distributions with $D_1$ being involutive, $D_1 \not\subset \mathcal{C}(D_2)$, and a vector field $f$ such that $[f, D_0]\subset D_1$ and $D_2 = D_1 + [f, D_1]$. Given two vector fields $v_1,v_2$ such that $D_1=D_0+\Span{v_1,v_2}$, a necessary condition for the existence of a vector field $v_c \in D_1, v_c\notin D_0$ such that the distributions $H_1=D_0+\Span{v_c}$ and $H_2 = D_1 + \Span{[f,v_c]}$ satisfy $H_1\subset\Cauchy{H_2}$ is that
    \begin{equation}\label{eq5:lem_cond}
        (\alpha^1)^2[v_1,[v_1,f]] + 2\alpha^1\alpha^2[v_1,[v_2,f]]+(\alpha^2)^2[v_2,[v_2,f]] \in D_2
    \end{equation}
    admits a nontrivial solution $\alpha^1,\alpha^2$. Candidates for $v_c$ are then obtained as $v_c=\alpha^1v_1+\alpha^2v_2 \text{ mod }D_0$.
\end{lem}
\smallskip

For the detailed proof, we refer to \cite{GstottnerAnalysisControlFlat2023}, where it is also shown that \eqref{eq5:lem_cond} admits at most two independent non-trivial solutions. Consequently, there exist at most two independent vector fields of the form $v_c = \alpha^1 v_1 + \alpha^2 v_2$ -- which are not collinear modulo $D_0$\footnote{Two vector fields $v,w$ are called collinear modulo a distribution $D$, if there exist functions $\alpha,\beta$, not both identically zero, such that $\alpha v + \beta w \in D$.} -- that meet \eqref{eq5:lem_cond}. 
Next, we show that $v_c$ of Lemma \ref{lem_1} plays a key role in the refined algorithm introduced later. Formulated in the adapted coordinates corresponding to the GTF, the following lemma highlights this role by linking the vector field $v_c$ from Lemma~\ref{lem_1} to the unit vector fields associated with subsystem $\Xi_3$ of \eqref{eq4:triang_form_1}.

\begin{lem}
    \label{lem_2}
    Consider a system of the form~(\ref{eq4:triang_form_1}) with vector fields~(\ref{eq4:vec_fields_a_b}). Define $ D_1 = \Span{b_1, b_2} $ and $D_{i+1}=D_i+[a,D_i]$, $i=1,\ldots,p$, and assume that $D_1,\ldots,D_p$ are involutive and $D_{p+1}$ is non-involutive. Further, let \mbox{$\underset{\inv}{D_{p-1}} \crktwo \underset{\inv}{D_{p}} \crktwo \underset{\notinv}{D_{p+1}}$}, satisfy the assumptions of Lemma~\ref{lem_1}. Then, a suitable choice for the vector field $ v_c $ from Lemma~\ref{lem_1} is given by $ v_c = \pad{z^{n - (p - 1)}} $.
\end{lem}
\begin{IEEEproof}
    Since the vector field $ f $ from Lemma~\ref{lem_1} corresponds to the drift vector field $ a $ of \eqref{eq4:triang_form_1}, a possible choice for the vector fields $v_1,v_2$ from Lemma~\ref{lem_1} is given by $ v_1 = \ad[p-1]{a}{b_1} $ and $ v_2 = \pad{z^{n - (p - 1)}} $, respectively. Substituting $v_2$ in the third term of (\ref{eq5:lem_cond}) yields
    \begin{equation*}
        [\pad{z^{n - (p - 1)}}, [\pad{z^{n - (p - 1)}}, a]]  \in D_{p+1},
    \end{equation*}
    which directly follows from the triangular structure of $a$ as indicated in \eqref{eq4:vec_fields_a_b}.
\end{IEEEproof}

We now present a refined procedure for constructing a sequence of feedback-invariant distributions (\ref{eq:alg_1_sequence}).

\begin{alg}\label{alg_new}
    A sequence \eqref{eq:alg_1_sequence} is computed as follows:
    \begin{enumerate}[label=\Alph*)]
        \item If $D_i$ is involutive, then $D_{i+1}=D_i+[f,D_i]$.
        \item If $D_i$ is non-involutive and $\Cauchy{D_i} = D_{i-1}$, then \mbox{$D_{i+1}=\bar{D}_i$}.
        \item If $D_i$ is non-involutive, $\Cauchy{D_i} \neq D_{i-1}$, and
        \begin{equation}\label{eq5:dist_case_C}
            \left(\Cauchy{D_i}\cap D_{i-2} \right)\crktwo \underset{\inv}{D_{i-1}} \crktwo \underset{\notinv}{D_i}
        \end{equation}
        satisfy the assumptions of Lemma~\ref{lem_1}, then proceed with:
        \begin{enumerate}[label=\roman*)]
            \item If a vector field $v_c$ with $\alpha_1,\alpha_2$ as a nontrivial solution of \eqref{eq5:lem_cond} exists, then replace $D_i$ with $D_{i} = D_{i-1} + \Span{[f, v_c]}$. If $D_i$ is involutive, continue with $D_{i+1}=D_i+[f,D_i]$, otherwise continue with $D_{i+1}=\bar{D}_i$.
            \item If no vector field $v_c$ exists, then $D_{i+1}=\bar{D}_i$.
        \end{enumerate}
        \item If $D_i$ is non-involutive, $\Cauchy{D_i} \neq D_{i-1}$, and there exist no distributions of the form \eqref{eq5:dist_case_C} satisfying the assumptions of Lemma~\ref{lem_1}, then $D_{i+1}=\bar{D}_i$.
    \end{enumerate}
\end{alg}

Basically, Algorithm~\ref{alg_new} extends Algorithm~\ref{alg_1} by introducing the update step C) i), which replaces the set of directions $[f,D_{i-1}]$ by selectively adding only $\Span{[f,v_c]}$. This refinement is advantageous in cases with small relative degrees, since it delays the inclusion of directions $\pad{z^{k_1}}$ or $\pad{z^{k_1+k_2}}$, see the discussed after \eqref{eq5:vec_fields_b1}. Therefore, it significantly broadens the class of systems for which $x$-flat output candidates can be identified, e.g., to systems SFE to \eqref{eq4:triang_form_2} with $q\leq 1$, to systems that are static feedback linearizable after a one/two-fold prolongation and beyond. Nevertheless, Algorithm~\ref{alg_new} does not provide a complete solution for identifying $x$-flat outputs of systems SFE to the GTF in the general case.

\section{Examples}\label{sec6}

This section presents examples that illustrate our main contributions. First, we show that an $x$-flat system \eqref{eq1:ai_sys} which is not SFE to the GTF can be rendered SFE by at most $n-2$ prolongations, as stated in Theorem~\ref{thm:prolong}. In the presented example, however, a one-fold prolongation already suffices. 
\begin{exa}
    Consider the system
    \vspace{-2ex}\\\vspace{1ex}
    \begin{minipage}[t]{0.25\textwidth}
    \vspace{2ex}
    \begin{align*}
        \dot{x}^1 & = u^1, \\
        \dot{x}^2 & = x^3 + x^4 u^1,
    \end{align*}
    \end{minipage}
    \hspace{-2.5em}
    \begin{minipage}[t]{0.2\textwidth}
    \begin{align*}
        \dot{x}^3 & = \alpha + (\beta - x^5) u^1, \\
        \dot{x}^4 & = x^5 + \gamma u^1, \\
        \dot{x}^5 & = u^2 \;
    \end{align*}
    \end{minipage}\\
    from \cite{NicolauNormalFormsXflat2022}, where the functions \( \alpha, \beta, \gamma \) may depend on \( x^1, \ldots, x^4 \). The system admits the flat output \( y = (x^1, x^2) \) with relative degree \( K = (1,1) \) and differential difference \( d = 3 \). As shown in \cite{GstottnerTriangularNormalForm2024}, it is not SFE to the GTF, since not all of the corresponding codistributions \eqref{eq4:seq_codist} are integrable. However, the system obtained by prolonging each input once is SFE to the GTF, because the codistributions
    \begin{equation*}
        \begin{aligned}
            Q_K = Q_{(1,1)} & = \Span{\D x^1, \D x^2, \D u^1, \D x^3 + u^1 \D x^4}, \\
            Q_{K+1} = Q_{(2,2)} & = \Span{\D x^1, \D x^2, \D x^3, \D x^4, \D u^1}, \\
            Q_{K+2} = Q_{(3,3)} & = \Span{\D x^1, \D x^2, \D x^3, \D x^4, \D x^5, \D u^1}, \\
            Q_{K+3} = Q_{(4,4)} & = \Span{\D x_e}
        \end{aligned}
    \end{equation*}
    associated with the prolonged system, whose state is given by $x_e = (x^1, \ldots, x^5, u^1, u^2)$, are integrable.
\end{exa}
We proceed with a well-known example, namely the planar VTOL aircraft, as considered, e.g., in \cite{FliessLieBacklundApproachEquivalence1999}.
\begin{exa}
    Consider the model of the planar VTOL aircraft
    \begin{equation}\label{eq6:vtol}
        \arraycolsep=3pt
        \begin{array}{rclcrcl}
            \dot{x} & = & v_{x}, & & \dot{v}_{x} & = & \epsilon\cos(\theta)u^2 - \sin(\theta)u^1, \\ [0.2cm]
            \dot{z} & = & v_{z}, & & \dot{v}_{z} & = & \cos(\theta)u^1 + \epsilon\sin(\theta)u^2-1, \\ [0.2cm]
            \dot{\theta} & = & \omega, & & \dot{\omega} & = & u^2,
        \end{array}
    \end{equation}
    with the drift vector field $f = v_x\pad{x} +  v_z\pad{z} + \omega\pad{\theta} - \pad{v_{z}}$ and the input vector fields $g_1 = -\sin(\theta)\pad{v_{x}} + \cos(\theta)\pad{v_{z}}$ and $g_2 = \epsilon\cos(\theta)\pad{v_{x}} + \epsilon\sin(\theta)\pad{v_{z}} + \pad{\omega}$.
    Applying Algorithm~\ref{alg_1} yields the sequence of distributions
    \begin{equation*}
        \underset{\inv}{D_1} \crktwo \underset{\notinv}{D_2} \crktwo \underset{\inv}{D_3} = \TX,
    \end{equation*}
    with $\Cauchy{D_2}=\emptyset$, providing no useful information for identifying flat outputs. In contrast, applying Algorithm~\ref{alg_new} yields a sequence \eqref{eq:alg_1_sequence} that enables the identification of flat-output components. Since $D_1$ is involutive, we begin with case (A): \\[1ex] 
    \noindent\textbf{Step 1 (A):} We compute $D_2 = D_1 + [f, D_1]$. Since $D_2$ is non-involutive, $\Cauchy{D_2} \neq D_1$, and the distributions \eqref{eq5:dist_case_C} given by \mbox{$\emptyset \crktwo D_1 \crktwo D_2 $} satisfy the assumptions of Lemma \ref{lem_1}, we continue with (C). \\[1ex]
    \noindent\textbf{Step 2 (C):} Condition \eqref{eq5:lem_cond} reduces to $\alpha^1\alpha^2 = 0$. Thus, with $\alpha^1=1,\alpha^2=0$ and $\alpha^1=0,\alpha^2=1$, we obtain two candidates: $v_{c,1}=g_1$ and $v_{c,2}=g_2$. This means, we have a branch in our algorithm and we have to continue with two possible distributions $D_{2,i}=D_1+\Span{[f,v_{c,i}]}$, $i=1,2$, that replace $D_2$ from Step 1. Both are non-involutive. Let $D_{3,i}=D^{(1)}_{2,i}=\bar{D}_{2,i}$ be their involutive closures. We proceed with case (A). \\[1ex]
    \noindent\textbf{Step 3 (A):} For both distributions, $D_{4,i}=D_{3,i}+[f,D_{3,i}]=\TX$ holds. Finally, we obtain two possible sequences
    \begin{equation*}
        \underset{\inv}{D_1} \underset{1}{\subset} \underset{\notinv}{D_{2,i}} \underset{1}{\subset} \underset{\inv}{D_{3,i}} \crktwo \TX, \quad i=1,2,
    \end{equation*}
    with the annihilators $D^{\perp}_{3,1} = \Span{\D \theta, \cot(\theta)\D x + \D z}$ and $D^{\perp}_{3,2} = \Span{\D x - \epsilon \cos(\theta)\D \theta, \D z - \epsilon \sin(\theta)\D \theta}$. Let us first consider the possible pair of functions $(\theta, x\cot(\theta) + z )$ corresponding to $D^{\perp}_{3,1}$ with its relative degree $K=(2,2)$. Assuming that $(\theta, x\cot(\theta) + z )$ is a valid $x$-flat output $\varphi$ of \eqref{eq6:vtol}, we obtain $R=(4,4)$ using \eqref{eq3:two_inputs_R_K_n}. However, it can be verified that \eqref{eq5:nec_suf_flat} is not satisfied, which implies that $\varphi$ is not a flat output.\footnote{Nevertheless, as shown in \cite{SchoberlCalculatingFlatOutputs2011}, the component $\theta$ belongs to a flat output whose second component depends on the input $u$.} In contrast, the functions $(x-\epsilon\cos(\theta), z+\epsilon\sin(\theta))$ corresponding to $D^{\perp}_{3,2}$ form a valid $x$-flat output $\varphi$, see also, e.g., \cite{MartinDifferentLookOutput1996}, where this result was derived by physical considerations. This flat output is not only compatible with the GTF, but also with the form \eqref{eq4:triang_form_2} with $q=1$, as shown in \cite{GstottnerStructurallyFlatTriangular2022}. By means of the geometric characterization provided therein, the same flat output was derived.
\end{exa}

The subsequent example demonstrates the capability of our enhanced algorithm for the identification of $x$-flat output candidates, even for systems that are SFE to GTF, but not SFE to a triangular form with a known comprehensive characterization, such as \eqref{eq4:triang_form_2} with $q=0$ or $q=1$.


\begin{exa}\label{exa_3}
Consider the following system, which for illustrative purposes is already given in the form \eqref{eq4:triang_form_1}:
\begin{equation}\label{eq:exa_722}
    \arraycolsep = 1.4pt
    \begin{array}{rclcrclcrcl}
       \dot{z}^1 & = & z^2, & \hspace{1em} & \dot{z}^3 & = & z^4, & \hspace{1em} & \zdot^6 & = & z^7, \\
       \dot{z}^2 & = & v^1, & &  \dot{z}^4 & = & z^5(1+v^1), & & \zdot^7 & = & v^2. \\
        & & & &  \dot{z}^5 & = & z^6 + z^2v^1, & & \\
   \end{array}
\end{equation}
System \eqref{eq:exa_722} admits the $x$-flat output $\varphi=(z^1,z^3)$ with the relative degree $K=(2,2)$. Note that the flat output $\varphi$ is evident by inspection. However, in the absence of a complete geometric characterization, determining $\varphi$ in general coordinates is exceedingly complex. The drift vector field and the input vector fields are given by
\begin{equation}\label{eq:vec_fields_722}
   \begin{aligned}
       a & = z^2\pad{z^1} + z^4\pad{z^3} + z^5\pad{z^4} + z^6\pad{z^5} + z^7\pad{z^6} , \\
       b_1 & = \pad{z^2} + z^5\pad{z^4} + z^2\pad{z^5} , \quad b_2 = \pad{z^7} \; . \\
   \end{aligned}
\end{equation} 
Since the application of Algorithm~\ref{alg_1} for~\eqref{eq:exa_722} does not yield valid candidates for $x$-flat output components, we proceed by applying Algorithm~\ref{alg_new}. Given the fact that $D_1=\Span{b_1,b_2}$ is involutive, we again start with case (A): \vspace{1ex} \\
\noindent\textbf{Step 1 (A):}
We compute
\begin{equation*}
    D_2 = \Span{b_1,\, \partial_{z^7},\, \partial_{z^1} + z^5\,\partial_{z^3} + (z^2 - z^6)\,\partial_{z^4},\, \partial_{z^6}},
\end{equation*}
which is non-involutive. Its Cauchy characteristic distribution is \( \Cauchy{D_2} = \Span{\partial_{z^7}} \neq D_1 \). Since the assumptions of Lemma~\ref{lem_1} are satisfied for the sequence \( \emptyset \crktwo D_1 \crktwo D_2 \), we proceed with case~(C). \\[-1ex]

\noindent\textbf{Step 2 (C):}  
With the vector fields $ f = a $, $ v_1 = b_1 $, and $ v_2 = b_2 $ as defined in \eqref{eq:vec_fields_722}, condition~\eqref{eq5:lem_cond} reduces to $\alpha^1 (z^2\,\partial_{z^3} + \partial_{z^4}) \in D_2$. This yields the -- up to a multiplicative factor -- unique solution \( \alpha^1 = 0 \), \( \alpha^2 = 1 \), and hence \( v_c = v_2 = \partial_{z^7} \). Next, we replace $D_2$ from Step 1 with \mbox{$D_2 = D_1 + \Span{[a, v_c]} = \Span{b_1,\, \partial_{z^7},\, \partial_{z^6}}$}, which is involutive. Hence, we compute $D_3=D_2+[a,D_2]$, which is non-involutive. Since $\Cauchy{D_3} \neq D_2$ applies and the sequence \eqref{eq5:dist_case_C} takes the form $\Span{\partial_{z^7}} \crktwo D_2 \crktwo D_3$, satisfying the assumptions of Lemma~\ref{lem_1}, we proceed with case~(C). \\[-1ex]

\noindent\textbf{Step 3 (C):}  
With the vector fields $f=a$, $v_1 = \partial_{z^2} + z^5\,\partial_{z^4} + z^2\,\partial_{z^5}$, and $v_2 = \partial_{z^6}$, condition~\eqref{eq5:lem_cond} takes the form
\begin{equation*}
    \begin{aligned}
        \alpha^1(z^2\,\partial_{z^3} + \partial_{z^4}) - \alpha^1\alpha^2\,\partial_{z^4} \in \\
        \Span{b_1,\, \partial_{z^7},\, \partial_{z^6},\, \partial_{z^5},\, z^5\,\partial_{z^1} - (z^2 - z^6)\,\partial_{z^2} + (z^5)^2\,\partial_{z^3}}.
    \end{aligned}
\end{equation*}
Up to a scalar factor, the only solution of the condition given above is $\alpha^1 = 0$, $\alpha^2 = 1$, yielding $v_c = v_2 = \partial_{z^6}$. $D_3$ from Step 3 is then replaced with $D_3 = D_2 + [a, v_c] = \Span{b_1,\, \partial_{z^7},\, \partial_{z^6},\, \partial_{z^5}}$, with its involutive closure taking the form $D_4 = D_3^{(1)} = \bar{D}_3 = \Span{\partial_{z^7},\, \partial_{z^6},\, \partial_{z^5},\, \partial_{z^4},\, \partial_{z^2}}$. \\[-1ex]

\noindent\textbf{Step 4 (A):}  
Finally, we compute $D_5 = D_4 + [a, D_4] = \TX$, completing the construction of the distribution sequence
\begin{equation*}
    \underset{\inv}{D_1} \underset{1}{\subset} \underset{\inv}{D_2} \underset{1}{\subset} \underset{\notinv}{D_3} \underset{1}{\subset} \underset{\inv}{D_4} \crktwo D_5 = \TX.
\end{equation*}
Given that \( D_4^{\perp} = \Span{\D z^1,\, \D z^3} \), flat output candidates are functions of the form \( \varphi^i = \psi^i(z^1, z^3) \). Indeed, the pair \( \varphi = (z^1, z^3) \) defines a flat output for the system \eqref{eq:exa_722}.
\end{exa}

\section{Conclusion}

In this work, we showed that the GTF subsumes several known normal forms and proved that any $x$-flat system can be transformed into the GTF via finite input prolongations, as illustrated in Example 1. We also proposed a refined algorithm for identifying flat outputs, and demonstrated its applicability on both a physical system (Example 2) and an academic one (Example 3). While a full characterization of \eqref{eq4:triang_form_1} appears intractable, future work could focus on characterizing \eqref{eq4:triang_form_2} for $q>1$ and determining conditions under which Algorithm~\ref{alg_new} applies to systems SFE to this form.

\appendices

\section{}\label{app:thm_triang_forms} 
Proof of Theorem~\ref{thm:trian_forms_1_2}: We demonstrate that a system (\ref{eq4:triang_form_2}) can be transformed into the GTF by reordering the states and applying invertible state transformations. First, we reorder the equations of subsystem $\Sigma_1$ and those of the left-hand side integrator chain of $\Sigma_3$:
\begin{equation}\label{proof_triang_1}
\arraycolsep = 1.4pt
\begin{array}{rclcrclcrcl}
    \dot{z}^1_{1,1} & = & z^2_{1,1}, & & \dot{z}^1_{3,1} & = & z^2_{3,1}, & & \dot{z}^1_{1,2} & = & z^2_{1,2}, \\
    & \vdots && & & \vdots & & && \vdots & \\
    \dot{z}^{n_{1,1}}_{1,1} & = & z^1_2, & & \dot{z}^{n_3-q}_{3,1} & = & v^1, & & \dot{z}^{n_{1,2}}_{1,2} & = & z^2_2 \; . \\
    \dot{z}^1_2 & = & z^1_{3,1}, &&&  \\
\end{array}
\end{equation}
Next, consider the second equation of $\Sigma_2$. Given the triangular structure of the system, we introduce the state transformation $\bar{z}^3_2 = z^3_2 z^1_{3,1} + a^2_2(z_1, z^1_2, z^2_2, z^3_2)$ with its inverse $z^3_2 = \hat{\psi}^3_2(z_1,z^1_2, z^2_2, \bar{z}^3_2, z^1_{3,1})$. Differentiating $\bar{z}^3_2$ with respect to time gives
\begin{equation}\label{proof_1_time_deriv}
\dot{\bar{z}}^3_2 = \alpha^3_2(z_1, z^1_2, z^2_2, \bar{z}^3_2, z^4_2, z^1_{3,1}) + \hat{\psi}^3_2 z^2_{3,1},
\end{equation}
with $\partial_{z^4_2}\alpha^3_2 \neq 0$. We now apply a further state transformation $\bar{z}^4_2 = \alpha^3_2 + \hat{\psi}^3_2z^2_{3,1}$, which is invertible with respect to $z^4_2$, and differentiate $\zb^4_2$ with respect to time. Hence, we construct an integrator chain by successively introducing state transformations and computing their time derivatives. To demonstrate that the dimension $n_2$ of $\Sigma_2$ and the length $n_3-q$ of the left-hand side integrator chain of $\Sigma_3$ affect neither SFE to the GTF nor the relative degrees \eqref{thm:relative_degrees}, we consider the following case distinction: \vspace{1ex} \break 
1) Case $n_2 \geq n_3-q+2$: The process outlined above is applied to the states of $\Sigma_2$ until $v^1$ occurs explicitly, i.e., until
\begin{equation*}
    \begin{aligned}
        \zbdot^{n_3-q+1}_2 & \hspace{-0.3em} = \hspace{-0.2em}\alpha^{n_3-q+1}_2(\ldots, z^{n_3-q+2}_2, \ldots, z^{n_3-q-1}_{3,1}) + \hat{\psi}^3_2z^{n_3-q}_{3,1} \\
        & \hspace{-0.3em} = \zb^{n_3-q+2}_2, \\
        \dot{\bar{z}}^{n_3-q+2}_2 & \hspace{-0.3em} = \hspace{-0.2em}\alpha^{n_3-q+2}_2(\ldots, z^{n_3-q+3}_2, \ldots,z^{n_3-q}_{3,1}) + \hat{\psi}^3_2v^1,
    \end{aligned}
\end{equation*}
with $\partial_{z^{n_3-q+3}_2}\alpha^{n_3-q+2}_2 \neq 0$. 
Accordingly, the resulting equations are of the form\footnote{Note that for $n_2=n_3-q+2$, the last equation of~\eqref{proof_triang_2} is of the form $\dot{\zb}^{n_2}_2=\alpha^{n_2}_2(z_1, \bar{z}_2,z_{3,1},z^1_{3,2}) + \hat{\psi}^3_2v^1$, with $\partial_{z^1_{3,2}}\alpha^{n_2}_2\neq0$.}
\begin{equation}\label{proof_triang_2}
    \begin{array}{rcl}
        \zdot^2_2 & = & \zb^3_2, \\
        & \vdots &  \\
         \dot{\bar{z}}^{n_3-q+2}_2 & = & \hspace{-0.2em}\alpha^{n_3-q+2}_2(\ldots, z^{n_3-q+3}_2, \ldots,z^{n_3-q}_{3,1}) + \hat{\psi}^3_2v^1, \\
         \zdot^{n_3-q+3}_2 & = & z^{n_3-q+4}_2z^1_{3,1} + \bar{a}^{n_3-q+3}_2(\ldots, z^{n_3-q+4}_2), \\
         & \vdots & \\
         \zdot^{n_2}_2 & = & z^1_{3,2} + \bar{B}(z_1, \zb_2)z^1_{3,1}, \\
    \end{array}
\end{equation}
with $\zb_2=(z^1_2, z^2_2, \zb^3_2, \ldots, \zb^{n_3-q+2}, z^{n_3-q+3}, \ldots, z^{n_2}_2)$. Consequently, by appending (\ref{proof_triang_2}) and the right-hand side integrator chain of $\Sigma_3$ to (\ref{proof_triang_1}), we obtain a system in GTF. \vspace{1ex} \break 
2) Case $n_2 < n_3-q+2$: We proceed by successively introducing state transformations and computing their first-order time derivatives, until we reach
\begin{equation*}
    \dot{\bar{z}}^{n_2}_2 = \alpha^{n_2}_2(z_1, \zb_2, z^1_{3,1}, \ldots, z^{n_2-2}_{3,1}, z^1_{3,2} ) + \hat{\psi}^3_2 z^{n_2-1}_{3,1}
\end{equation*}
where $\zb_2=(z^1_2, z^2_2, \zb^3_2, \ldots, \zb^{n_2}_2)$ and $\partial_{z^1_{3,2}}\alpha^{n_2}_2\neq 0 $. Next, we introduce the state transformation $\bar{z}^1_{3,2} = \alpha^{n_2}_2 + \hat{\psi}^3_2 z^{n_2-1}_{3,1}$. Repeating the procedure of introducing state transformations and computing their time derivatives, we continue for the first $n_3-q-n_2+2$ states of $z_{3,2}$, until we reach
\begin{equation*}
    \begin{aligned}
        \dot{\bar{z}}^{n_3-q-n_2+2}_{3,2} & =  \alpha^{n_3-q-n_2+2}_{3,2} + \hat{\psi}^3_2 v^1,
    \end{aligned}
\end{equation*}  
where $  \alpha^{n_3-q-n_2+2}_{3,2}$ may depend on $z_1, \zb_2, z_{3,1}, \zb^1_{3,2},\ldots,$ $\zb^{n_3-q-n_2+2}_{3,2}$ and satisfies $\partial_{z^{ n_3-q-n_2+3 }_{3,2}}\alpha^{n_3-q-n_2+2}_{3,2} \neq 0$. As a result, we obtain
\begin{equation}\label{proof_triang_3}
    \arraycolsep=1.4pt
    \begin{array}{rclcrcl}
        && && \zbdot^1_{3,2} & = & \zb^2_{3,2}, \\
         \zdot^2_2 & = & \zb^3_2, & \hspace{1em} & & \vdots & \\
        \dot{\bar{z}}^3_2 & = & \bar{z}^4_2, && \zbdot^{n_3-q-n_2+1}_{3,2} & = & \zb^{n_3-q-n_2+2}_{3,2}, \\
        & \vdots & && \zbdot^{n_3-q-n_2+2}_{3,2} & = & \alpha^{n_3-q-n_2+2}_{3,2} + \hat{\psi}^3_2 v^1, \\
        \zbdot^{n_2}_2 & = & \zb^1_{3,2}, && \dot{z}^{n_3-q-n_2+3}_{3,2} & = & z^{n_3-q-n_2+4}_{3,2}, \\
        && && & \vdots & \\
        &&&& \dot{z}^{n_3}_{3,2} & = & v^2. \\
    \end{array}
\end{equation}
By appending \eqref{proof_triang_3} to \eqref{proof_triang_1}, the overall system is given in GTF. Finally,  analyzing the lengths of the integrator chains in both cases 1) and 2), we conclude that the relative degrees are given by \eqref{thm:relative_degrees}.

\section{}\label{app:prolong} 
Proof of Theorem~\ref{thm:prolong}: Consider an $x$-flat system to which the static feedback transformation (\ref{eq3:input_trf}) has been applied, followed by a $d$-fold input prolongation, yielding the prolonged system \vspace{-1ex}
\begin{equation}\label{eq:prf_prolong_system}
    \arraycolsep = 1.4pt
    \begin{array}{rcl}
         \dot{x} & = & \bar{f}(x) + \bar{g}_1(x)\bar{u}^1 + \bar{g}_2(x)\bar{u}^2, \\
         \dot{\bar{u}}^1_{[0,d-1]} & = & \bar{u}^1_{[1,d]}, \quad \dot{\bar{u}}^2_{[0,d-1]} = \bar{u}^2_{[1,d]} \; ,
    \end{array}
\end{equation}
with the new input $\bar{u}_{[d]}$. Obviously, the relative degrees and the multi-index from Definition \ref{def:flatness} are now given by $K+d$ and $R+d$, respectively. To prove that (\ref{eq:prf_prolong_system}) is SFE to the GTF, it is necessary to verify that the codistributions $Q_{R-1}, \ldots, Q_{R+d-1}$, cf. (\ref{eq4:seq_codist}), are all integrable. The derivatives of the flat output up to the order $R+d-1$ follow the structure of (\ref{eq3:output_deriv}), extended by:
\begin{equation*}
    \arraycolsep = 1.4pt
    \begin{array}{rclcrcl}
        \varphi^1_{[r_1+1]} & = & \bar{u}^1_{[d+1]}, & & \varphi^2_{[r_2+1]} & = &\varphi^2(x,\bar{u}^1_{[0,d+1]}, \bar{u}^2_{[0,1]}), \\
        & \vdots & & & & \vdots & \\ 
        \varphi^1_{[r_1+d-1]} & = & \bar{u}^1_{[2d-1]}, & &\varphi^2_{[r_2+d-1]} & = &\varphi^2(x,\bar{u}^1_{[0,2d-1]}, \bar{u}^2_{[0,d-1]}). \\
    \end{array}
\end{equation*}
Given that flatness of the original system implies $\Span{\D x} \subset \Span{\D \varphi_{[0,R-1]}} $ and $\Dim{\Span{\D \varphi_{[0,R-1]}}} = r_1+r_2= n + d$, it follows that 
\begin{equation*}
    \begin{aligned} 
        \Span{\D \varphi_{[0,R-1]}} & = \Span{\D x, \D \bar{u}^1_{[0,d-1]}}, \\
        \Span{\D \varphi_{[0,R]}} &= \Span{\D x, \D \bar{u}^1_{[0,d]}, \D \bar{u}^2}, \\
        & \hspace{0.5em} \vdots \\
        \Span{\D \varphi_{[0,R+d-1]}} & = \Span{\D x, \D \bar{u}^1_{[0,2d-1]}, \D \bar{u}^2_{[0,d-1]}} \; . \\
    \end{aligned}
\end{equation*}
Consequently, the associated codistributions $Q_{R-1}, \ldots, Q_{R+d-1}$ are indeed all integrable.

\bibliography{IEEEabrv, mybibfile}

\begin{thebibliography}{10}
\providecommand{\url}[1]{#1}
\csname url@samestyle\endcsname
\providecommand{\newblock}{\relax}
\providecommand{\bibinfo}[2]{#2}
\providecommand{\BIBentrySTDinterwordspacing}{\spaceskip=0pt\relax}
\providecommand{\BIBentryALTinterwordstretchfactor}{4}
\providecommand{\BIBentryALTinterwordspacing}{\spaceskip=\fontdimen2\font plus
\BIBentryALTinterwordstretchfactor\fontdimen3\font minus \fontdimen4\font\relax}
\providecommand{\BIBforeignlanguage}[2]{{%
\expandafter\ifx\csname l@#1\endcsname\relax
\typeout{** WARNING: IEEEtran.bst: No hyphenation pattern has been}%
\typeout{** loaded for the language `#1'. Using the pattern for}%
\typeout{** the default language instead.}%
\else
\language=\csname l@#1\endcsname
\fi
#2}}
\providecommand{\BIBdecl}{\relax}
\BIBdecl

\bibitem{FliessSurLesSystemes1992}
M.~Fliess, J.~L{\'e}vine, P.~Martin, and P.~Rouchon, ``Sur les syst{\`e}mes non lin{\'e}aires diff{\'e}rentiellement plats,'' \emph{C.R. Acad. Sci. Paris}, vol. 315, pp. 619--624, 1992.

\bibitem{FliessFlatnessDefectNonlinear1995}
------, ``Flatness and defect of non-linear systems: Introductory theory and examples,'' \emph{International Journal of Control}, vol.~61, no.~6, pp. 1327--1361, 1995.

\bibitem{FliessLieBacklundApproachEquivalence1999}
M.~Fliess, J.~Levine, P.~Martin, and P.~Rouchon, ``A {{Lie-B{\"a}cklund}} approach to equivalence and flatness of nonlinear systems,'' \emph{IEEE Transactions on Automatic Control}, vol.~44, no.~5, pp. 922--937, 1999.

\bibitem{GstottnerTrackingControlFlat2024}
C.~Gst{\"o}ttner, B.~Kolar, and M.~Sch{\"o}berl, ``Tracking control for (x, u)-flat systems by quasi-static feedback of classical states,'' \emph{Symmetry, Integrability and Geometry: Methods and Applications (SIGMA)}, vol.~20, no.~71, 2024.

\bibitem{MartinFeedbackLinearizationDriftless1994}
P.~Martin and P.~Rouchon, ``Feedback linearization and driftless systems,'' \emph{Mathematics of Control, Signals, and Systems}, vol.~7, no.~3, pp. 235--254, 1994.

\bibitem{NicolauFlatnessMultiInputControlAffine2017}
F.~Nicolau and W.~Respondek, ``Flatness of multi-input control-affine systems linearizable via one-fold prolongation,'' \emph{SIAM Journal on Control and Optimization}, vol.~55, no.~5, pp. 3171--3203, 2017.

\bibitem{NicolauFlatnessTwoinputControlaffine2016}
------, ``Flatness of two-input control-affine systems linearizable via a two-fold prolongation,'' in \emph{2016 {{IEEE}} 55th {{Conference}} on {{Decision}} and {{Control}} ({{CDC}})}.\hskip 1em plus 0.5em minus 0.4em\relax Las Vegas, NV, USA: IEEE, 2016, pp. 3862--3867.

\bibitem{GstottnerFiniteTestLinearizability2021}
C.~Gst{\"o}ttner, B.~Kolar, and M.~Sch{\"o}berl, ``A finite test for the linearizability of two-input systems by a two-dimensional endogenous dynamic feedback,'' in \emph{2021 {{European Control Conference}} ({{ECC}})}, 2021, pp. 970--977.

\bibitem{GstottnerNecessarySufficientConditions2023}
------, ``Necessary and sufficient conditions for the linearisability of two-input systems by a two-dimensional endogenous dynamic feedback,'' \emph{International Journal of Control}, vol.~96, no.~3, pp. 800--821, 2023.

\bibitem{BououdenTriangularCanonicalForm2011}
S.~Bououden, D.~Boutat, G.~Zheng, J.-P. Barbot, and F.~Kratz, ``A triangular canonical form for a class of 0-flat nonlinear systems,'' \emph{International Journal of Control}, vol.~84, no.~2, pp. 261--269, 2011.

\bibitem{GstottnerFlatTriangularForm2021}
C.~Gst{\"o}ttner, B.~Kolar, and M.~Sch{\"o}berl, ``On a flat triangular form based on the extended chained form,'' \emph{IFAC-PapersOnLine}, vol.~54, no.~9, pp. 245--252, 2021.

\bibitem{GstottnerStructurallyFlatTriangular2022}
------, ``A structurally flat triangular form based on the extended chained form,'' \emph{International Journal of Control}, vol.~95, no.~5, pp. 1144--1163, 2022.

\bibitem{SchoberlImplicitTriangularDecomposition2014}
M.~Sch{\"o}berl and K.~Schlacher, ``On an implicit triangular decomposition of nonlinear control systems that are 1-flat - a constructive approach,'' \emph{Automatica}, vol.~50, no.~6, pp. 1649--1655, 2014.

\bibitem{SilveiraFlatTriangularForm2015}
H.~B. Silveira, P.~S. Pereira Da~Silva, and P.~Rouchon, ``A flat triangular form for nonlinear systems with two inputs: {{Necessary}} and sufficient conditions,'' \emph{European Journal of Control}, vol.~22, pp. 17--22, 2015.

\bibitem{GstottnerTriangularNormalForm2024}
C.~Gst{\"o}ttner, B.~Kolar, and M.~Sch{\"o}berl, ``A triangular normal form for x-flat control-affine two-input systems,'' in \emph{2024 28th {{International Conference}} on {{Methods}} and {{Models}} in {{Automation}} and {{Robotics}} ({{MMAR}})}, 2024, pp. 298--303.

\bibitem{KolarPropertiesFlatSystems2016}
B.~Kolar, M.~Sch{\"o}berl, and K.~Schlacher, ``Properties of flat systems with regard to the parameterization of the system variables by the flat output,'' \emph{IFAC-PapersOnLine}, vol.~49, no.~18, pp. 814--819, 2016.

\bibitem{GstottnerLinearizationFlatTwoInput2020}
C.~Gst{\"o}ttner, B.~Kolar, and M.~Sch{\"o}berl, ``On the linearization of flat two-input systems by prolongations and applications to control design,'' \emph{IFAC-PapersOnLine}, vol.~53, no.~2, pp. 5479--5486, 2020.

\bibitem{GstottnerFlatSystemPossessing2023}
------, ``A flat system possessing no (x, u)-flat output,'' \emph{IEEE Control Systems Letters}, vol.~7, pp. 1033--1038, 2023.

\bibitem{KaiExtendedChainedForms2006}
T.~Kai, ``Extended chained forms and their application to nonholonomic kinematic systems with affine constraints: {{Control}} of a coin on a rotating table,'' in \emph{Proceedings of the 45th {{IEEE Conference}} on {{Decision}} and {{Control}}}, 2006, pp. 6104--6109.

\bibitem{LiCharacterizationFlatnessExtended2013}
S.~Li, C.~Xu, H.~Su, and J.~Chu, ``Characterization and flatness of the extended chained system,'' in \emph{Proceedings of the 32nd {{Chinese Control Conference}}}, 2013.

\bibitem{LiMultiinputControlaffineSystems2016}
S.~Li, F.~Nicolau, and W.~Respondek, ``Multi-input control-affine systems static feedback equivalent to a triangular form and their flatness,'' \emph{International Journal of Control}, vol.~89, no.~1, pp. 1--24, 2016.

\bibitem{TilburyGoursatNormalForms1994}
D.~Tilbury and S.~Sastry, ``On {{Goursat}} normal forms, prolongations, and control systems,'' in \emph{Proceedings of 1994 33rd {{IEEE Conference}} on {{Decision}} and {{Control}}}, vol.~2, 1994, pp. 1797--1802.

\bibitem{TilburyTrajectoryGenerationNtrailer1995}
D.~Tilbury, R.~Murray, and S.~Shankar~Sastry, ``Trajectory generation for the {{N-trailer}} problem using {{Goursat}} normal form,'' \emph{IEEE Transactions on Automatic Control}, vol.~40, no.~5, pp. 802--819, 1995.

\bibitem{MurraySteeringNonholonomicSystems1991}
R.~Murray and S.~Sastry, ``Steering nonholonomic systems in chained form,'' in \emph{{{Proceedings}} of the 30th {{IEEE Conference}} on {{Decision}} and {{Control}}}, 1991, pp. 1121--1126.

\bibitem{MurrayNonholonomicMotionPlanning1993}
------, ``Nonholonomic motion planning: Steering using sinusoids,'' \emph{IEEE Transactions on Automatic Control}, vol.~38, no.~5, pp. 700--716, 1993.

\bibitem{LiFlatOutputsTwoinput2012}
S.-J. Li and W.~Respondek, ``Flat outputs of two-input driftless control systems,'' \emph{ESAIM: Control, Optimisation and Calculus of Variations}, vol.~18, no.~3, pp. 774--798, 2012.

\bibitem{KolarContributionsDifferentialGeometric2017}
B.~Kolar, \emph{Contributions to the Differential Geometric Analysis and Control of Flat Systems}, ser. Modellierung Und {{Regelung}} Komplexer Dynamischer {{Systeme}}.\hskip 1em plus 0.5em minus 0.4em\relax Aachen: Shaker, 2017, vol.~35.

\bibitem{GstottnerAnalysisControlFlat2023}
C.~Gst{\"o}ttner, \emph{Analysis and Control of Flat Systems by Geometric Methods}, ser. {{Modellierung und Regelung komplexer dynamischer Systeme}}.\hskip 1em plus 0.5em minus 0.4em\relax D{\"u}ren: Shaker Verlag, 2023, vol.~59.

\bibitem{NicolauNormalFormsXflat2022}
F.~Nicolau, C.~Gst{\"o}ttner, and W.~Respondek, ``Normal forms for x-flat two-input control-affine systems in dimension five,'' \emph{IFAC-PapersOnLine}, vol.~55, no.~30, pp. 394--399, 2022.

\bibitem{SchoberlCalculatingFlatOutputs2011}
M.~Sch{\"o}berl and K.~Schlacher, ``On calculating flat outputs for {{Pfaffian}} systems by a reduction procedure - {{Demonstrated}} by means of the {{VTOL}} example,'' in \emph{2011 9th {{IEEE International Conference}} on {{Control}} and {{Automation}} ({{ICCA}})}, 2011, pp. 477--482.

\bibitem{MartinDifferentLookOutput1996}
P.~Martin, S.~Devasia, and B.~Paden, ``A different look at output tracking: Control of a vtol aircraft,'' \emph{Automatica}, vol.~32, no.~1, pp. 101--107, 1996.

\end{thebibliography}

\end{document}